\documentclass[12pt,a4paper]{article}

\usepackage[utf8]{inputenc}

\usepackage[english]{babel}

\usepackage{graphicx}
\usepackage{placeins} 
\usepackage{subfig}
\usepackage{wrapfig}

\usepackage{amsmath}
\usepackage{amsthm}
\usepackage{amsfonts}

\theoremstyle{plain}
\newtheorem{theorem}{Theorem}

\newtheorem{proposition}[theorem]{Proposition}

\theoremstyle{definition}

\newtheorem{remark}[theorem]{Remark}

\newcommand{\dd}{{\rm d}}

\newcommand{\la}{\langle}
\newcommand{\ra}{\rangle}
\renewcommand{\H}{\mathcal{H}}

\DeclareMathOperator\RES{Res}

\usepackage[style=alphabetic-verb,
backref=false,
hyperref=auto,
maxnames=4,
backend=biber]{biblatex}

\usepackage[babel]{csquotes}
\usepackage[colorlinks=true,
        linkcolor=black,
        citecolor=black,
        filecolor=black,
        urlcolor=black,
        bookmarks=true,
        bookmarksopen=false,
        plainpages=false,
        pdfpagemode=UseNone,
        pdfpagelabels=true]{hyperref}

\usepackage[color=cyan]{todonotes}
\begin{document}

\title{Rigorous a-posteriori analysis using numerical
    eigenvalue bounds in a surface growth model}
\author{Dirk Blömker, Christian Nolde}

\maketitle

\begin{abstract}
In order to prove numerically the global existence
and uniqueness of smooth solutions of a
fourth order, nonlinear PDE,
we derive rigorous a-posteriori upper bounds on the
supremum of the numerical range of the linearized operator.
These bounds also have to be
easily computable in order to be applicable to our rigorous a-posteriori
methods, as we use them in each time-step of the numerical discretization.
The final goal is to establish global
bounds on smooth local solutions, which then establish global uniqueness.
\end{abstract}

\section{Introduction}

This paper deals with the rigorous numerical verification of 
global existence and uniqueness
of smooth solutions to the surface growth equation
\begin{align}
u_t = - u_{xxxx} - \big( (u_x)^2 \big)_{xx}. \label{eqn:sfg}
\end{align}
on $x \in [0,2\pi]$ with periodic boundary conditions.

This equation, usually with additional lower order terms and noise,
was introduced as a phenomenological model for
the growth of vapor deposited amorphous surfaces
\cite{Si-Pl:94b,Ra-Li-Ha:00a},
and was also used
to describe ion-sputtering processes, where a surface is eroded
by an ion-beam \cite{Cu-Va-Ga:05}.
The one dimensional equation appears as a model for the boundaries of
terraces in the epitaxy of silicon \cite{frischverga06}.
A more detailed list of references can be found in the review 
article \cite{blomkerromito15}.

Analytically, this PDE was studied by Blömker and Romito
in several papers which are reviewed in \cite{blomkerromito15},
including the existence of smooth local solutions in the
largest critical space and an example for a blowup in the case
of the complex valued equation, which rules out the possibility 
that standard energy estimates
alone might be sufficient to proof global uniqueness.

Except for small initial data,
there are no analytic methods to prove the existence of smooth
global solutions known so far. The equation only has uniform in time bounds
on the spatial $L^2$-norm of solutions, and global existence of 
solutions for all
initial conditions in $L^2$. But in contrast to that uniqueness only holds for
initial conditions of higher regularity like $C^0$, $H^\frac{1}{2}$,
or some suitable Besov-space (see \cite{blomkerromito15} 
for details).

For problems where analytic methods are not able to produce results yet,
the application of rigorous computational methods is a steadily increasing
field over the recent years. The used methods vary as much as the problems
they are applied to.
For proving numerically the existence of
solutions for PDEs, in addition to our approach, there are methods based
on topological arguments
like the Conley index, see
\cite{paapeetal08, daylessardmischaikov07}, for example.
For solutions of elliptic PDEs there are methods using Brouwer's fixed-point
theorem, as discussed in the review article \cite{Pl:08} and
the references therein. Finite element methods to 
obtain lower bounds on eigenvalues can be found in e.g.
\cite{rannacherheuveline01}, \cite{carstensengedicke14} and 
\cite{carstensengedicke14-2}.
For periodic solutions or invariant manifolds
for dissipative PDE see for example \cite{zgliczynski10} and 
\cite{bergmirelesreinhardt16}.
A nice introductory overview is \cite{berglessard15}.

\subsection{The previous worst case method}

Our method is based on \cite{robinsonetal07} 
which is formulated for the3D Navier-Stokes equation,
and related ideas can be found in \cite{morosipizzocchero08},
although no numerical experiments are present in these papers.
A different approach to the problem is studied by 
\cite{lessardcyranka17}, which is more in the direction of 
the methods cited in the previous section.

The key idea of \cite{robinsonetal07} is to establish 
a scalar ODE that bounds 
the difference $d$ between a unique smooth
local solution $u$
and an arbitrary approximation $\varphi$, which is provided 
by a numerical method, for instance.
As the existence and uniqueness of $u$ for the surface growth 
equation, for example in $H^{1}$, is 
known as cited before, we obtain the following result:
As long as we can bound the $H^1$-norm   $\|d_x\|$, we are able to
use the unique continuation of the smooth local solution and obtain a
unique smooth solution up to the blow up time of 
our error bound on $\|d_x\|$.

As for any initial value $u_0\in H^1$ there is a time 
$T^*(u_0)$ with the property
that if there was no blow up until time $T^*$,
there can not occur one afterwards
(see \cite{blomkernolderobinson15} Theorem 3, Time Condition).
Thus one can also obtain global existence and uniqueness
by controlling the $H^1$-norm of the error up to that time.
Similar properties are also well known for 3D Navier-Stokes.

Let us comment in more detail on the result of 
\cite{blomkernolderobinson15}, which is closely related 
to \cite{robinsonetal07}.
The key analytic result of that paper is the following 
differential inequality for the error
\begin{align}
\label{e:compode}
    \partial_t \| d_x \|^2 \leq
    \frac{7^7}{2} \| d_x \|^{10} +
    \Big( 18 \| \varphi_{xx} \|_\infty^2 - \frac{1}{2}
    \Big) \| d_x \|^2
    + 2 \| \RES \|_{-1}^2,
\end{align}
where $\RES := \varphi_t + \varphi_{xxxx} + ({\varphi_x}^2)_{xx}$ is 
the residual of the approximation $\varphi$ that measures
how close $\varphi$ is to being a solution of \eqref{eqn:sfg}.

As the coefficients
of the right hand side of \eqref{e:compode} depend only on the numerical data,
using the time discretization of the numerical solution $\varphi$ 
this ODE could be
evaluated rigorously for instance by using interval arithmetic.
As we were mainly interested in performing a case study whether 
the approach is working at all,
we did not yet implement interval arithmetic in our numeric simulations,
but this is just a technical issue in programming.

Further, in \cite{blomkernolderobinson15} we showed that this 
approach could give global existence
for initial conditions larger than the analytic smallness result, which 
is limited to solutions of $H^1$-norm smaller than $1/2$.
In the numerical simulations we could easily treat larger initial 
conditions like $u_0 = \sin(x)$.
On the other hand, the method based on \eqref{e:compode} still 
fails for even moderately increased
frequencies in the initial value (without dampening by the amplitude) 
like $u_0 = \sin(2x)+ \cos(3x)$,
as the $H^1$-norm gets too large.

Let us finally remark that due to the scaling properties of the equation,
we can always treat some initial conditions that are arbitrarily 
large in $H^1$.
If $u(t,x)$ is any spatially $2\pi$-periodic solution of \eqref{eqn:sfg},
then for any $k \in \mathbb{N}$ the rescaled solution 
$u_k(t,x) = u(k^4t,kx)$ is also a $2\pi$-periodic solution. But now it 
is easy to see that for the initial condition
$\|u_k(0,\cdot)\|_{H^1}\to \infty$ if $k\to\infty$.

\subsection{Improvement based on numerical eigenvalues}

It turned out in our numerical experiments of \cite{blomkernolderobinson15},
that the most sensitive part for our rigorous method 
based on \eqref{e:compode} is
the $18 \| \varphi_{xx} \|_\infty^2$ term that leads to a 
strong exponential growth.
In contrast to that the residual $\RES$ seems to be always extremely small,
indicating a fast convergence of the numerical method we use 
to obtain $\varphi$. But we are analytically far from proving any 
convergence of the numerical method.

As the quintic nonlinearity in our ODE for the error 
\eqref{e:compode} immediately 
leads to a blow up in finite time, once the error is sufficiently large,
we were looking for a way to improve our error estimate, by 
replace our previous ``worst case'' estimates leading to the 
term $18 \| \varphi_{xx} \|_\infty^2$.
This estimate was purely analytic and largely relied on general 
interpolation inequalities,
bounding the respective quadratic form of the linearized operator.
Hereby, we are following the idea of
\cite{nakaohashimoto09, nakaokinoshitakinura12}, where
the spectrum of the linearized operator
 is analyzed. In our case this is the non-symmetric
\[L_\varphi u =  -\partial_x^4 u - 2\partial_x^2 ( \varphi_x  u_x )\;,
\]
where $\varphi$ is some given numerical data , and thus $L_\varphi u$
is just the linearization of the full nonlinear SPDE 
(\ref{eqn:sfg}) along the numerical approximation $\varphi$.

The bound is based on a rigorous numerical method 
for the largest eigenvalue, 
which in the case of an
unstable linear operator yields substantially better results,
at the price of a significantly higher computational time.

Let us comment in more detail on this. In order to derive an 
improvement of \eqref{e:compode},
we are interested in the supremum of the numerical range 
of $L_\varphi$, which means
we want to bound the quadratic form
\begin{equation*}
\lambda(\varphi) = \sup_{\|u_x\| = 1 }\langle \partial_x L_\varphi u,
\partial_x  u \rangle
\end{equation*}
in order to finally obtain a bound
\[
\langle \partial_x L_\varphi u, \partial_x  u \rangle
\leq \lambda(\varphi) \|u_x\|^2.
\]
This is equivalent to bounding the largest eigenvalue of the 
symmetrized operator $\frac12(L_\varphi+L_\varphi^\star)$.
Although there are already results for upper bounds on the largest
eigenvalue of self-adjoint operators (for example \cite{liu15}), 
we have the requirement that our estimate
is also (relatively) easy and fast to compute in order to be applicable to our
a-posteriori method as it has to be calculated in every time step
of the discretization.
\subsection{Structure of the paper}

In Section \ref{sec:setting}, we state the basic notation 
used throughout the paper.
The main result for the numerical eigenvalue is stated 
in Section \ref{sec:MThm},
and proven in Section \ref{sec:proof}.
In Section \ref{sec:comp} we compare the new estimate 
with the previous worst case estimate
and demonstrate how much better the verification for global 
existence and uniqueness works with the new estimate based 
on the numerical eigenvalue.

\section{Setting \& Problem} \label{sec:setting}
As solutions to our surface growth equation (\ref{eqn:sfg}) are
subject to periodic boundary conditions on $[0,2\pi]$ with mean average
zero, we are working on the Hilbert space
\[
\mathcal{H}= \Big \{ u : \mathbb{R} \to \mathbb{R} \ :\
2\pi \text{-periodic},\  \int_0^{2\pi}u(x)\; \; \dd x=0 \Big\}
\]
with standard $L^2$-scalar product $\la \cdot,\cdot \ra$ 
and corresponding $L^2$-norm
\[
\|u\|  = \Big(\int_0^{2\pi}|u(x)|^2 \; \; \dd x \Big)^{1/2}\;.
\]
We further define the Sobolev-spaces
\[
\mathcal{H}^k= \{u \in \mathcal{H}\ :\ \partial_x^k u \in L^2([0,2\pi]) \}.
\]
Note that by periodicity $u \in \mathcal{H}^1$ implies  $u_x \in \mathcal{H}$.
Moreover, we have Poincare-inequality with optimal constant $1$
\[
\|u\| \leq \| u_x\|  \quad\text{for all }u\in \mathcal{H}^1
\]
and thus $\|u_x\|$ is a norm on $\mathcal{H}^1$, equivalent
to the standard $H^1$-Sobolev norm.

Furthermore, interpolation inequality holds also with constant $1$
\[
\|u_x\|^2 \leq \| u_{xx}\|\|u\|  \quad\text{for all }u\in \mathcal{H}^2\;.
\]
In both cases the constants are easy to compute. 
For details see \cite{nolde17}.

Let us recall in more detail the results of \cite{blomkernolderobinson15}.
There, in order to control the $\mathcal{H}^1$-norm of a unique 
smooth local solution $u$ to the surface growth equation (\ref{eqn:sfg}),
we derived a differential inequality to bound the
$\mathcal{H}^1$-norm of the difference
\[
d(x,t) := u(x,t) - \varphi(x,t),
\]
where
$\varphi$ is any arbitrary, but sufficiently smooth approximation, that
satisfies periodic boundary conditions.
In the numerical examples we always use a
spectral Galerkin method in space and a semi-implicit Euler 
scheme in time, which we then extend by piece-wise linear interpolation
of the numerical data in time. Thus $\varphi$ is arbitrarily 
smooth in space (i.e., $C^\infty$) and Lipschitz (i.e., $W^{1,\infty}$) in time.

Using a standard a-priori type estimate, the
differential inequality for the error is given by
\begin{align}
\tfrac12 \partial_t \| d_x \|^2 &=
    \underbrace{ \la d_{xx} , d_{xxxx} + 2 (d_x \varphi_x)_{xx} \ra }
        _{\text{A}+\text{B}}
    + \underbrace{ \la d_{xx} , ({d_x}^2)_{xx} \ra }
        _{\text{C}}
    + \underbrace{ \la d_{xx}, \RES \ra }
        _{\text{D}} \label{eqn:sfgODE}\\
&\leq \tfrac{7^7}{4} \| d_x \|^{10}
    + \big( 9 \| \varphi_{xx} \|_{L^\infty}^2 - \tfrac14 \big) \| d_x \|^2
    + \| \RES \|_{H^{-1}}^2 \nonumber,
\end{align}
with residual $\RES:= \varphi_t + \varphi_{xxxx} + ({\varphi_x}^2)_{xx}.$
The estimate above is based on a crude ``worst case'' estimate 
for A+B  and was established in \cite{blomkernolderobinson15}.

Our aim of this paper is to improve this estimate specifically for the term A+B,
by using  a numerical calculation
that computes a more problem specific estimate.

Therefore consider the linearized operator
\[
L_\varphi u =  -\partial_x^4 u - 2\partial_x^2 ( \varphi_x  u_x ).
\]
We are interested in bounding the quadratic form
\begin{equation}
\lambda = \sup_{\|u_x\| = 1}\langle \partial_x L_\varphi u,
\partial_x  u \rangle
\label{eqn:quadForm}
\end{equation}
in order to finally obtain a bound
\[
\text{A}+\text{B} =
\langle \partial_x L_\varphi d, \partial_x  d \rangle
\leq \lambda\|d_x\|^2.
\]
Note that we neglect the explicit dependence of $\lambda$ on 
$\varphi$ and thus on time in the notation.

In order to transform this to an eigenvalue problem in $L^2$, we
 substitute $v=u_x$ in \eqref {eqn:quadForm} and immediately get
\[
\lambda =  \sup_{\|v\| = 1}\langle A_\varphi v, v \rangle  \;.
\]
with non-symmetric  operator
\begin{align}
A_\varphi u =  -\partial_x^4 u - 2\partial_x^3 ( \varphi_x  u ).
\label{def:operatorA}
\end{align}
For the numerical computation of $\lambda$ we also use a 
spectral Galerkin method.
Define $H_n$ as the $2n$-dimensional subspace spanned by $e^{ix}$,
\ldots, $e^{i n x}$ and its complex conjugates $e^{-ix}$,
\ldots, $e^{-i n x}$. Note that we can omit the constant 
mode due to our solution space
$\mathcal{H}$.
Denote by $P_n$ the orthogonal projection onto $H_n$.

Finally, we set the numerical approximation of $\lambda$ as
\begin{equation}
\lambda_n :=  
\sup_{\|u\| = 1, u\in\mathcal{H} }\langle P_n A_\varphi P_n u, u \rangle 
=  
\sup_{\|u\| = 1, u\in H^n }\langle A_\varphi u, u \rangle
\end{equation}
which is just the largest eigenvalue of a symmetric
$2n \times 2n$ matrix given by the symmetrized matrix 
$\frac12(P_n A_\varphi P_n+ P_nA^* P_n)$.

Obviously, as the supremum is over a larger set, it immediately holds that
\[
\lambda_n \leq \lambda.
\]
and moreover, $\lambda_n$ is monotone and thus convergent.

In the following sections we want to bound $\lambda$ from above 
by $\lambda_n$ plus an explicit error term,
which is the difficult task.
%

\section{Main Theorem}
\label{sec:MThm}
First, let us recall the ``worst case'' estimate from
\cite{blomkernolderobinson15}.
\begin{proposition} \label{prop:worstcase}
    Consider $A_\varphi$ as defined in (\ref{def:operatorA}) with
    $\varphi \in W^{2,\infty}$, then it holds that
    \[
    \langle A_\varphi u, u \rangle
    \leq -\tfrac12\|u_{xx}\|^2  + \tfrac92 \| \varphi_{xx}\|^2_\infty \|u\|^2
    \leq [-\tfrac12  + \tfrac92 \| \varphi_{xx}\|^2_\infty ] \cdot\|u\|^2
    \]
    for all $u \in H^2$.
\end{proposition}
Note that we are working with smooth local solutions or finite Fourier series,
so this estimate will only be applied to sufficiently smooth $u$.
\begin{proof}
The estimate is first proven for sufficiently smooth 
$u \in H^4$, as the quadratic form needs a fourth derivative,
and then the estimate is easily extended by continuity of the 
quadratic form to $u \in H^2$.

First using integration by parts
\begin{align*}
\langle A_\varphi u, u \rangle
    &= -\|u_{xx}\|^2 + 2 \int \varphi_x u u_{xxx}\; \dd x  \\
&= -\|u_{xx}\|^2 - 2 \int \varphi_{xx} u u_{xx}\; \dd x
    - 2 \int \varphi_{x} u_x u_{xx}\; \dd x \\
&= -\|u_{xx}\|^2 - 2 \int \varphi_{xx} u u_{xx}\; \dd x
    + \int \varphi_{xx} u_{x}^2\; \dd x.  \\
\intertext{Now, H\"older, interpolation, and
    Poincare inequalities are used to obtain}
&\leq  -\|u_{xx}\|^2 + 2 \| \varphi_{xx}\|_\infty  \|u\|\|u_{xx}\|
    + \|\varphi_{xx}\|_\infty  \|u_{x}\|^2 \\
&\leq   -\|u_{xx}\|^2 + 3 \| \varphi_{xx} \|_\infty  \|u\|\|u_{xx}\| \\
&\leq  - \tfrac12 \|u_{xx}\|^2 + \tfrac92 \|\varphi_{xx}
    \|^2_\infty \|u\|^2 \\
&\leq -\tfrac12\|u\|^2  + \tfrac92 \| \varphi_{xx}\|^2_\infty \|u\|^2.
\end{align*}
\end{proof}
Thus we obtain for the supremum of the quadratic form defined in
(\ref{def:operatorA})
\begin{align}
\lambda \leq - \frac12  + \frac92\|\varphi_{xx}\|^2_\infty\;.
    \label{ineq:worstcase}
\end{align}
This is the worst case estimate used in 
\cite{blomkernolderobinson15} to obtain the differential 
inequality stated in \eqref{eqn:sfgODE}.

Instead, the following theorem shows an improved estimate
by analyzing the quadratic form (\ref{eqn:quadForm}) separately
for different mode ranges.
\begin{theorem} \label{thm:evbound}
    Let $u$ be a smooth local solution to our surface growth equation
    (\ref{eqn:sfg}) with initial condition 
    $u(0)\in\mathcal{H}^1$, $\varphi$ an arbitrary  $H_n$-valued
    approximation and $H_n$, $\lambda$ and $\lambda_n$ be defined as in
    Section \ref{sec:setting}.
    Then, for
    \[
    n \geq \sqrt2 C_\varphi
    = \sqrt2 ( 2 \|\varphi_{xxx}\|_\infty
    + 6  \| \varphi_{xx} \|_\infty
    + 4  \| \varphi_{x}  \|_\infty )
    \]
    it holds that
    \[
    \lambda_n \leq \lambda
    \leq \lambda_n + \frac12 \max \Big\{ 2 C_\varphi^2
    \frac{9\|\varphi_{xx}\|_\infty^2 -2\lambda_n  }{n^2}  \ ,\
    9\|\varphi_{xx}\|_\infty^2 + 2 \lambda_n- \frac12 n^4 \Big\}\;.
    \]
\end{theorem}
\begin{remark}
 Note that due to monotonicity $\lambda_n$ converges, and the previous 
 result shows the convergence of $\lambda_n$ to $\lambda$.
 Moreover, we obtain the asymptotic rate of convergence 
 \[
 \lambda = \lambda_n + \mathcal{O}(1/n^2)\;.
 \]
 On the other hand, for a given $n$ and a given $\varphi$, we can 
 calculate $\lambda_n$ and the error given by the previous 
 theorem fairly quickly.
\end{remark}

%
\section{Proof of the Theorem}
\label{sec:proof}
%
As a preparation, we split $u= p+q$, where $p \in H_n$ and $q \perp H_n$.
Thus
\begin{align*}
\lambda &= \sup_{\|u\| = 1 } \la A_\varphi u, u \ra \nonumber \\
& = \sup_{\|p\|^2+\|q\|^2 = 1 } \Big\{\la A_\varphi p, p \ra
+ \la A_\varphi p, q \ra + \la A_\varphi q, p \ra
+ \la A_\varphi q, q \ra   \Big\}.
\end{align*}
Now, we will treat these scalar products separately, where
we will denote with ``low modes'' the parts only depending on $p$
and with ``high modes'' everything solely depending on $q$.

Note that $A_\varphi$ is not symmetric and thus 
$\la A_\varphi p, q \ra \not= \la A_\varphi q, p \ra$, in general.
%
\subsection*{Low modes}
First, notice that by the brute force estimate of
Proposition \ref{prop:worstcase} we have
\[
\langle A_\varphi p, p \rangle
\leq -\frac12 \|p_{xx}\|^2 + \frac92 \|\varphi_{xx}\|^2_\infty\|p\|^2 \;.
\]
Second, it holds by definition of $\lambda_n$, as $p\in H_n$
\[
\langle A_\varphi p, p \rangle \leq \lambda_n \|p\|^2\;.
\]
In summary, we get for some $\eta_n \in[0,1]$, that we will fix later,
\begin{align*}
\langle A_\varphi p, p \rangle
\leq (1-\eta_n)\lambda_n \|p\|^2
- \frac12 \eta_n \|p_{xx}\|^2
+  \frac92 \eta_n \|\varphi_{xx}\|^2_\infty\|p\|^2.
\end{align*}
We do not use only the numerical eigenvalue to bound the quadratic form,
as we also need to control terms involving $\|p_{xx}\|$ arising 
in the estimate of  the mixed terms.
%
%
\subsection*{Mixed terms}
For the mixed terms we use the elementary estimates
\begin{align}
\|p\| \leq \|p_x\| \leq \|p_{xx}\|
\qquad\text{and}\qquad
\|q\| \leq \frac1n \|q_x\| \leq \frac1{n^2} \|q_{xx}\|
\label{def:impPoincare}
\end{align}
Note that any derivatives of $p$ and $q$ are still orthogonal in $\H$,
so the only terms in the mixed terms that are non-zero 
are the ones that contain $\varphi$.

We obtain first
\begin{align*}
\la A_\varphi p, q \ra
&= -2 \int ( \varphi_x  p )_{xxx}   q \; \dd x 
    =  2 \int ( \varphi_x  p )_{xx}   q_x \; \dd x \\
&= 2 \int ( \varphi_{xxx} p + 2\varphi_{xx} p_x + \varphi_x p_{xx} )
q_x \; \dd x \\
&\leq  2 \|q_x\| \cdot  \big( \|\varphi_{xxx}\|_\infty\|p\|
+ 2 \|\varphi_{xx}\|_\infty\|p_x\|
+ \|\varphi_{x}\|_\infty\|p_{xx}\| \big) \\
&\leq C_\varphi^{(1)} \frac1n \|q_{xx}\| \|p_{xx}\|
\end{align*}
with
\[C_\varphi^{(1)} = 2 \|\varphi_{xxx}\|_\infty +
4  \|\varphi_{xx}\|_\infty + 2 \|\varphi_{x}\|_\infty ].
\]
For the second mixed term we derive similarly
\begin{align*}
\la A_\varphi q , p  \ra
&= 2  \int ( \varphi_x  q )    p_{xxx} \; \dd x
= -2 \int ( \varphi_x  q )_x   p_{xx} \; \dd x \\
&\leq 2\|p_{xx}\| \cdot  \big( \|\varphi_{xx}\|_\infty\|q\|
+ \|\varphi_{x}\|_\infty\|q_x\|  \big) \\
&\leq C_\varphi^{(2)}  \frac1n \|q_{xx}\| \|p_{xx}\|
\end{align*}
with
\[C_\varphi^{(2)}= 2 \|\varphi_{xx}\|_\infty +
2 \|\varphi_{x}\|_\infty .
\]
Further, we define
\[
C_\varphi=C_\varphi^{(1)} + C_\varphi^{(2)}
= 2 \|\varphi_{xxx}\|_\infty + 6  \|\varphi_{xx}\|_\infty
+ 4 \|\varphi_{x} \|_\infty\;.
\]
%
%
\subsection*{High modes}
Finally, for the high modes we have no other option, but to 
use the rough ``worst case'' estimate
of Proposition \ref{prop:worstcase} which yields
\[
\langle A_\varphi q, q \rangle
\leq  - \frac12\|q_{xx}\|^2 + \frac92\|\varphi_{xx}\|^2_\infty \|q\|^2.
\]
We will apply the improved Poincare inequality (\ref{def:impPoincare})
which is valid on the high modes in a later step.
%
\subsection*{Summary}
%
Combining all estimates,  we obtain (using Young inequality
$ab\leq\frac12a^2+\frac12b^2$ and eliminating $p_{xx}$ terms)
\begin{align*}
\langle A_\varphi u, u \rangle
&=  \langle A_\varphi p, p \rangle
+ \langle A_\varphi p, q \rangle
+ \langle A_\varphi q, p \rangle
+ \langle A_\varphi q, q \rangle  \\
&\leq   (1-\eta_n)\lambda_n \|p\|^2
- \frac12 \eta_n \|p_{xx}\|^2
+  \frac92 \eta_n \|\varphi_{xx}\|^2_\infty \|p\|^2 \\
& \qquad + C_\varphi \frac1n \|q_{xx}\| \|p_{xx}\| \\
& \qquad - \frac12 \|q_{xx}\|^2
+ \frac92 \|\varphi_{xx}\|^2_\infty\|q\|^2 \\
&\leq  (1-\eta_n) \lambda_n \|p\|^2
+ \frac92 \eta_n \|\varphi_{xx}\|^2_\infty \|p\|^2 \\
& \qquad +\frac12 \Big(  \frac{C_\varphi^2}{n^2\eta_n} - 1\Big)
\|q_{xx}\|^2 +  \frac92 \|\varphi_{xx}\|_\infty^2\|q\|^2
\end{align*}
In order to apply the improved Poincare inequality (\ref{def:impPoincare})
for $q$, we define
\[
\eta_n := 2 \frac{C_\varphi^2}{n^2}
\; \text{and thus we need } n\geq \sqrt2 C_\varphi \; 
\text{to assert }\eta_n \leq 1.
\]
We obtain
\begin{align*}
\langle A_\varphi u, u \rangle & \leq  \Big[(1-\eta_n) \lambda_n
+ \frac92 \eta_n \|\varphi_{xx}\|^2_\infty \Big] \|p\|^2 \\
& + \frac12 \Big[9\|\varphi_{xx}\|^2_\infty
- \frac12 n^4 \Big]\|q\|^2
\end{align*}
which proves our main theorem
\begin{align*}
\lambda &= \sup_{\|u\| = 1 }\langle A_\varphi u, u \rangle
= \sup_{\|p\|^2+\|q\|^2 = 1 }\langle A_\varphi u, u \rangle \\
&\leq  \max \Big\{ \Big[(1-\eta_n) \lambda_n
+ \frac92 \eta_n \|\varphi_{xx}\|^2_\infty \Big]\ ,\
\frac12 \Big[9\|\varphi_{xx}\|_\infty^2 - \frac12 n^4
\Big]  \Big\} \\
&= \lambda_n + \frac12 \max
\Big\{\eta_n [ 9\|\varphi_{xx}\|^2_\infty
-2 \lambda_n  ]  \ ,\  9\|\varphi_{xx}\|^2_\infty
+ 2 \lambda_n- \frac12 n^4 \Big\}.
\end{align*}
\qed
%
\section{Simulations}
\label{sec:comp}
Before we come to the results of the simulations, 
let us first explain the numerical methods and necessary preparations
that we use to calculate 
$\varphi$ and the upper bounds on $\| d_x \|^2$.
\subsection*{Calculating $\varphi$}
To compute our arbitrary approximation $\varphi$, we use 
a spectral Galerkin method to convert the PDE to a system 
of ODEs.
Note that we only need any approximation, so no interval arithmetic is 
necessary in this step. 
The basis of eigenfunctions is in our case the standard Fourier 
basis $e_k = \frac{1}{\sqrt{2 \pi}}\exp(ikx)$.
As a welcome side effect this allows us  
to compute quantities like $L^2$ scalar products 
and norms very efficiently and accurately.

With $u := \sum_k a_k(t) e_k$ our surface growth equation (\ref{eqn:sfg})
turns into the following infinite system of coupled 
(through the nonlinearity) ODEs
\begin{align*}
a_k'(t) = - (ik)^4 a_k(t) - (ik)^2 
\underbrace{
    \Big( \sum_{s+l = k} (is) a_s(t) 
    \times (il)a_l(t) \Big)}_{b_k(t)} 
\quad \forall k.
\end{align*} 
For the spectral Galerkin approximation, we truncate the sum for $b_k$ 
to a finite range of modes.
To solve this system, 
we now use a semi-implicit Euler scheme, i.e. we use time 
$t_{j+1} = t_j + h$ in the 
linear part, and $t_j$ inside the nonlinearity (we could not solve for 
$t_{j+1}$)
\begin{align*}
\frac{1}{h} (a_k(t_{j+1}) 
- a_k(t_j) ) &=
- (ik)^4 a_k(t_{j+1}) 
- (ik)^2 b_k(t_j) \\
\intertext{and thus}
a_k(t_{j+1}) &=
(1 + h(k)^4)^{-1} (a_k(t_j) + h k^2 b_k(t_j) )
\end{align*}
for all $k$.
\subsection*{Applying the eigenvalue estimate}
Before we define how to calculate the bound on $\| d_x \|^2$, 
we have to incorporate the eigenvalue estimate from 
Theorem \ref{thm:evbound} into the bounding ODE (\ref{eqn:sfgODE}), 
which is given by 
\begin{align*}
\tfrac12 \partial_t \| d_x \|^2 &=
\underbrace{ \la d_{xx} , d_{xxxx} + 2 (d_x \varphi_x)_{xx} \ra }
_{\text{A}+\text{B}}
+ \underbrace{ \la d_{xx} , ({d_x}^2)_{xx} \ra }
_{\text{C}}
+ \underbrace{ \la d_{xx}, \RES \ra }
_{\text{D}} \\
&\leq \tfrac{7^7}{4} \| d_x \|^{10}
+ \big( 9 \| \varphi_{xx} \|_{L^\infty}^2 - \tfrac14 \big) \| d_x \|^2
+ \| \RES \|_{H^{-1}}^2.
\end{align*}
Let us denote the eigenvalue bound from Theorem \ref{thm:evbound}
with $\tilde{\lambda}$. If we want to apply this result 
to our framework, we have to consider, that 
in order to control the (C) and (D) terms, we need some part of the (A) term of
\begin{align*}
\frac{1}{2} \partial_t \|d_x\|^2 = &
\underbrace{    
    \la d_{xx}, d_{xxxx} + 2 (d_x \varphi_x)_{xx} \ra
}_{\text{A+B}}
+ 
\underbrace{ 
    \la d_{xx}, ({d_x}^2)_{xx} + \mathrm{Res} \ra 
}_{\text{C+D}}.
\end{align*}
Therefore, we split the first term into two parts ($\delta \in (0,1)$)
\begin{align*}
\frac{1}{2} \partial_t \|d_x\|^2 &=    
(1 - \delta) \la d_{xx}, d_{xxxx} + 2 (d_x \varphi_x)_{xx} \ra
+ \delta \la d_{xx}, d_{xxxx} + 2 (d_x \varphi_x)_{xx} \ra \\
&+ \la d_{xx}, ({d_x}^2)_{xx} + \mathrm{Res} \ra .  
\end{align*}
Now, we can bound the first term with our new method and 
the remaining parts like before in (\ref{eqn:sfgODE}).
If we do not fix the constants used in the Young inequalities, we have
\begin{align*}
\mathrm{A} &=- \| d_{xxx} \|^2 \\
|\mathrm{B}| & \leq   \varepsilon_B \| d_{xxx} \|^2
+ \frac{9}{4 \varepsilon_B} \| d_x \|^2 
\|	\varphi_{xx} \|_{\infty}^2 \\
|\mathrm{C}| &\leq  \varepsilon_C \| d_{xxx} \|^2
+ \frac{(\frac47 \varepsilon_C)^{-7}}{4} \| d_x \|^{10} \\
|\mathrm{D}| & \leq \varepsilon_D \| d_{xxx} \|^2 
+ \frac{1}{4 \varepsilon_D} \| \RES \|_{-1}^2,
\end{align*}
where we can set all $\varepsilon_{\{B,C,D\}} > 0$ arbitrary small.

In this case, our differential inequality is
\begin{align*}
\frac{1}{2} \partial_t \| d_x \|^2 
&\leq (1 - \delta) \tilde{\lambda} \| d_x \|^2 
+ \frac{9}{4 \varepsilon_B} \delta 
\| d_x \|^2	\| \varphi_{xx} \|_{\infty}^2
+ \frac{(\frac47 \varepsilon_C)^{-7}}{4} \| d_x \|^{10} \\
&+ \frac{1}{4 \varepsilon_D} \| \RES \|_{-1}^2 
+ \big( \delta \varepsilon_B + \varepsilon_C + \varepsilon_D - \delta \big) 
\| d_{xxx} \|^2,
\end{align*}
where $\varepsilon_{\{B,C,D\}} > 0$ and $ \delta \in (0,1)$. 
By substituting $\varepsilon_{\{C,D\}} := \delta \varepsilon_{\{C,D\}}$, 
this is equivalent to
\begin{align*}
\frac{1}{2} \partial_t \| d_x \|^2 
&\leq (1 - \delta) \tilde{\lambda} \| d_x \|^2 
+ \frac{9}{4 \varepsilon_B} \delta 
\| d_x \|^2	\| \varphi_{xx} \|_{\infty}^2
+ \frac{(\frac47 \delta \varepsilon_C)^{-7}}{4} \| d_x \|^{10} \\
&+ \frac{1}{4 \delta \varepsilon_D} \| \RES \|_{-1}^2 
+ \delta \big( \varepsilon_B + \varepsilon_C + \varepsilon_D - 1 \big) 
\| d_{xxx} \|^2,
\end{align*}
where $\varepsilon_{\{B,C,D\}} > 0$ and $ \delta \in (0,1)$. Next, 
we set $\varepsilon_B + \varepsilon_C + \varepsilon_D = 1$ to remove 
the last term, and therefore, our final ODE is given by
\begin{align}
\begin{split}
\frac{1}{2} \partial_t \| d_x \|^2 
&\leq (1 - \delta) \tilde{\lambda} \| d_x \|^2 
+ \frac{9 \delta}{4 \varepsilon_B}  
\| d_x \|^2	\| \varphi_{xx} \|_{\infty}^2
+ \frac{7^7}{4^8 (\delta \varepsilon_C)^7} \| d_x \|^{10}  \\
&+ \frac{1}{4 \delta \varepsilon_D} \| \RES \|_{-1}^2 
\label{eqn:newODE}
\end{split}
\end{align}
under the constraints $\varepsilon_{\{B,C,D\}} > 0$, 
$ \sum_{k \in \{B,C,D\}} \varepsilon_k = 1$, $\delta \in (0,1)$.
Unfortunately, there is no easy to determine global minimum 
in regard of the constraints.
We could rewrite this problem and finally solve it using 
Ferrari's method for quartic equations, but sadly
this approach has a very bad cost-benefit ratio as the involved calculations 
are too complex. 
Luckily, we can not do anything wrong here that breaks the rigorosity of our 
calculations, as valid parameter combinations just might not be optimal.
Therefore, we just use MATLAB's nonlinear optimization solver to find an 
approximate local minimum and update it after a given time interval 
(we could do this in every step, but given 
that the step-size is quite small and the data is continuous, 
this is not necessary and would just cost 
us lots of computational time) (see \cite{nolde17} for details).
\subsection*{Numerical Comparison}
We will now investigate the improvement of the new estimate
from Theorem \ref{thm:evbound} compared to the previous ``worst case''
estimate (\ref{ineq:worstcase}) in numerical simulations
of our rigorous a-posteriori method. Again, please note 
that interval arithmetic was not used for 
these simulations, and the results are therefore not rigorous.
We use the rigorous analytic bound for an ODE of the type 
\eqref{eqn:sfgODE} or \eqref{eqn:newODE} based on restarting the 
estimate on every time step.
Details of these calculations can be found in \cite{nolde17}.

\begin{figure}[p]%
    \centering
    \hspace*{\fill} %
    \subfloat[$u_0 =\sin(x)$, $N=128$]
    {\includegraphics[width=0.45\textwidth]
        {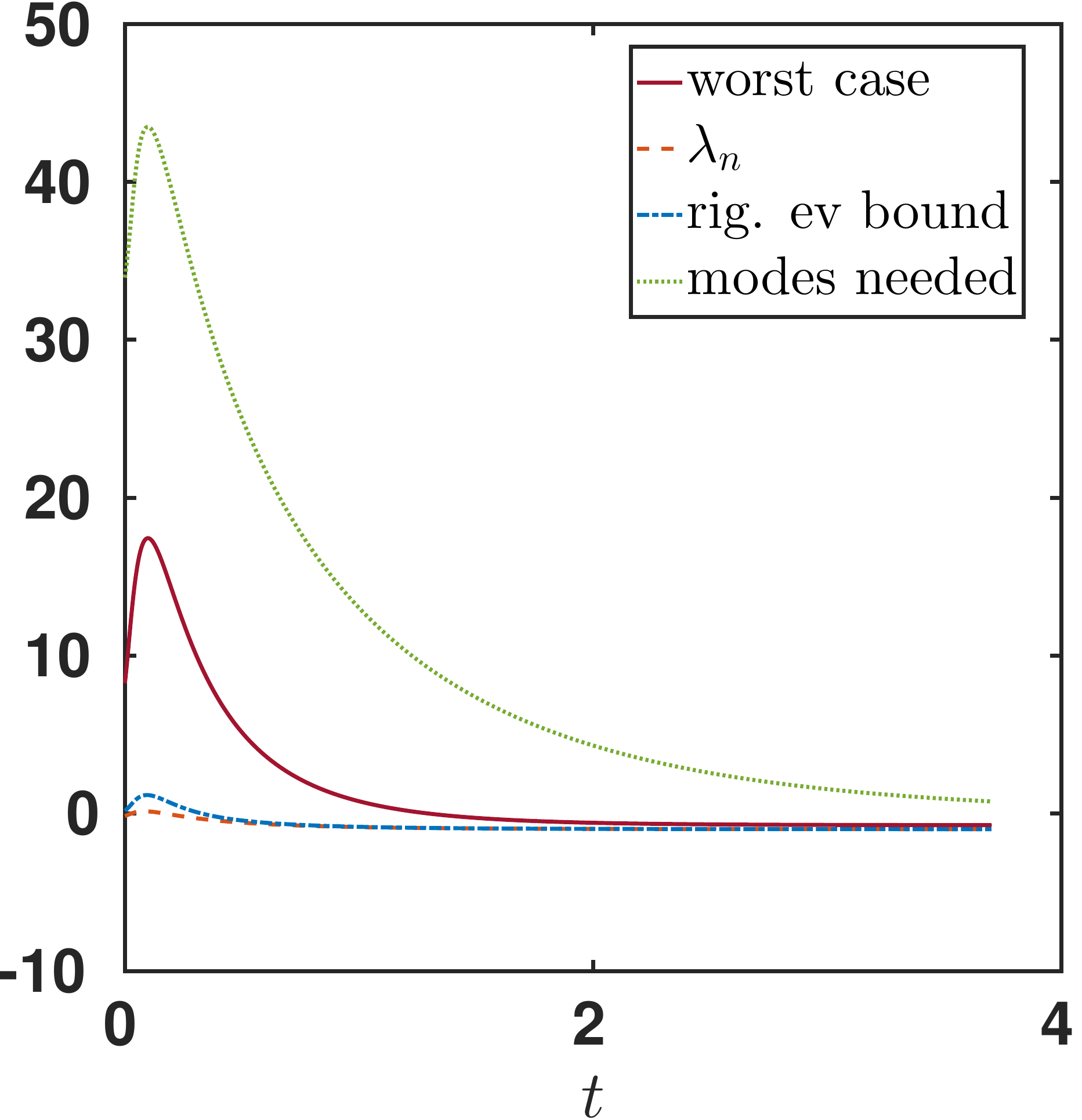}}
    \hspace*{\fill} %
    \subfloat[$u_0 = 2\sin(x)$, $N=256$]
    {\includegraphics[width=0.45\textwidth]
        {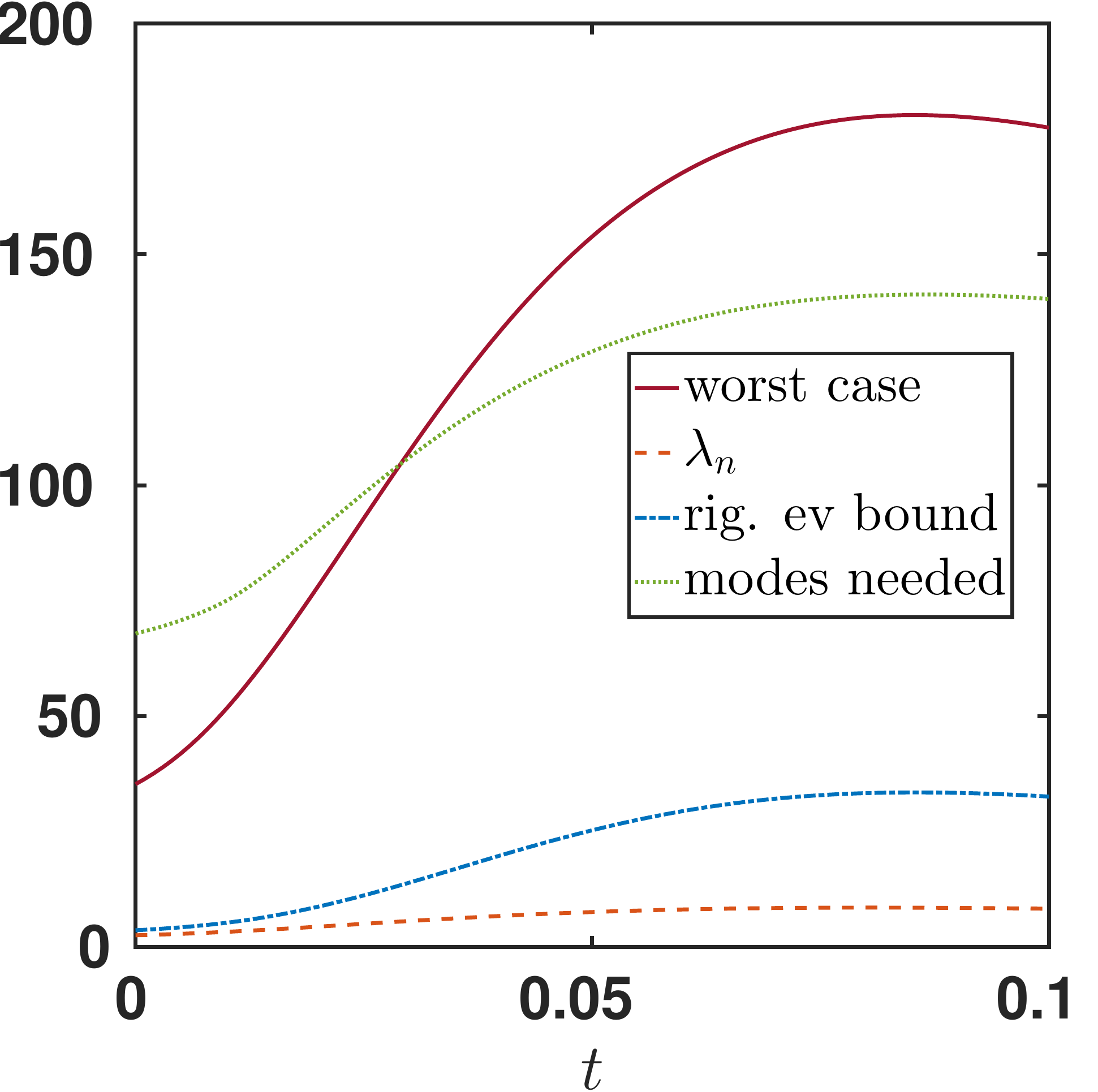} \label{img:evest2sin}}
    \hspace*{\fill} %
    \\
    \hspace*{\fill} %
    \subfloat[$u_0 = \sin(2x)+\cos(2x)$, $N=512$]
    {\includegraphics[width=0.45\textwidth]
        {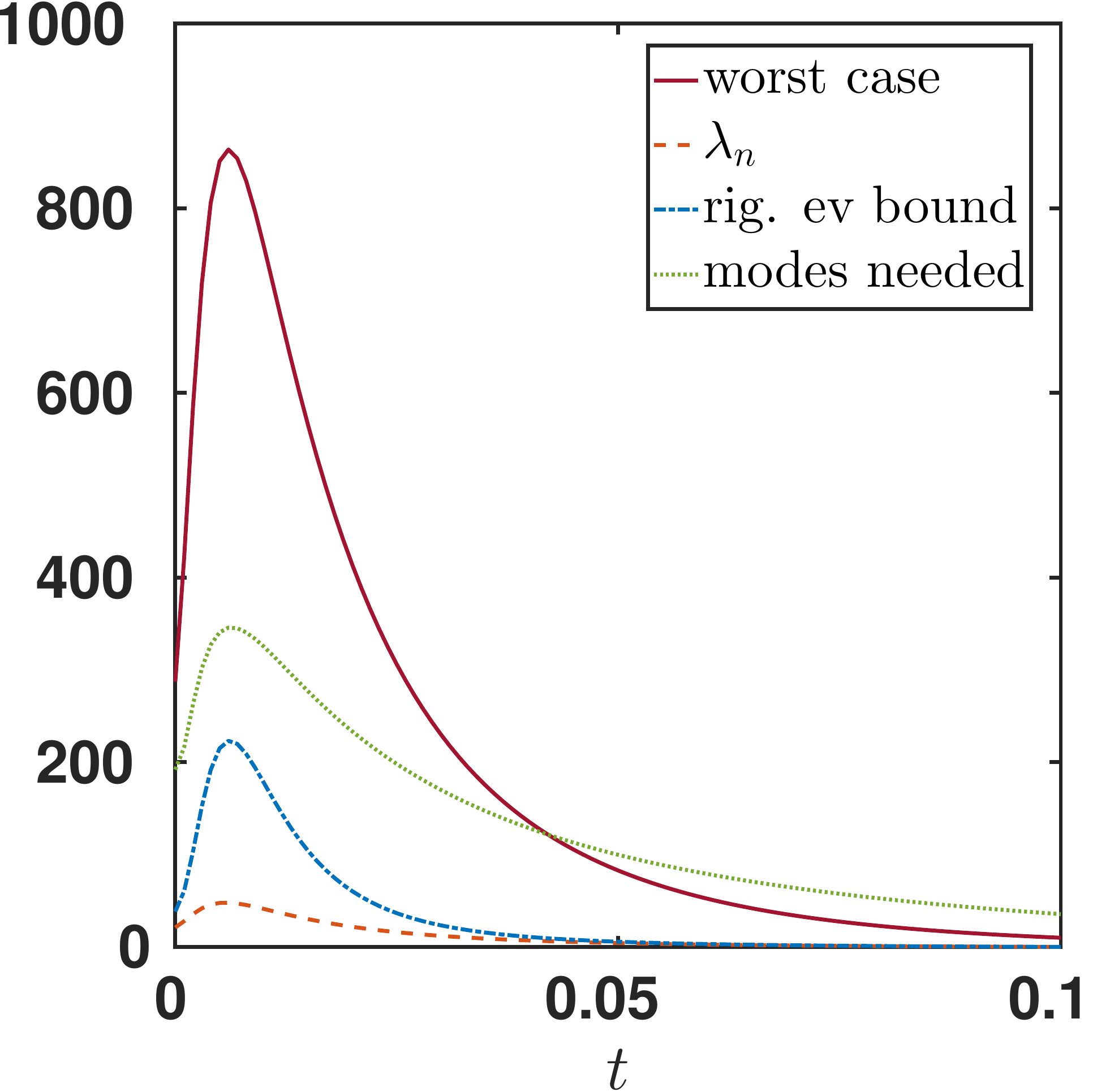}}
    \hspace*{\fill} %
    \subfloat[$u_0 = 1.5\sin(x)+\sin(2x)$, $N=512$]
    {\includegraphics[width=0.45\textwidth]
        {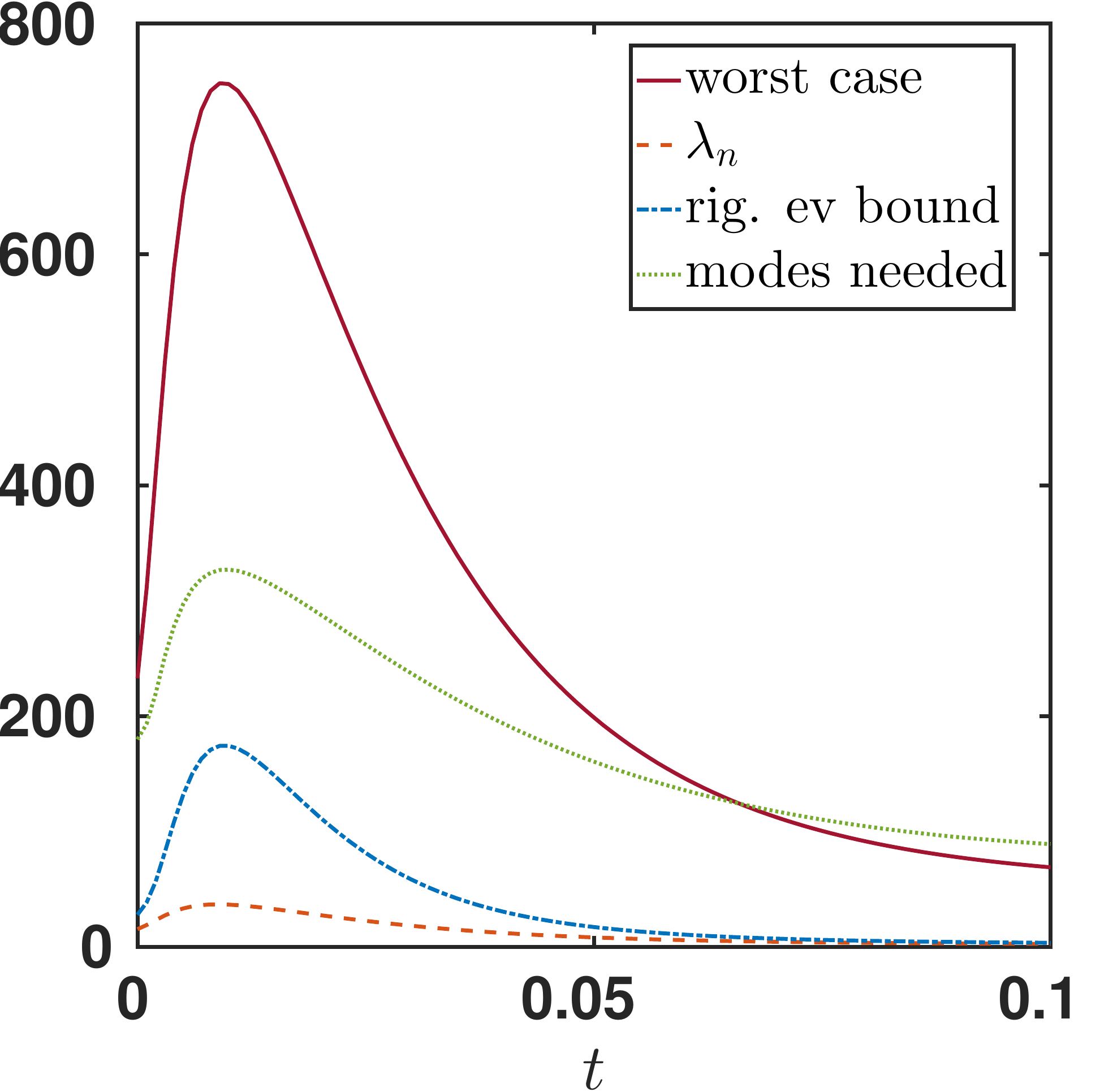}}
    \hspace*{\fill} %
    \\
    \caption
    {Comparison of the new Eigenvalue Estimate with the previous
    ``worst case''-estimate for different initial values.}
    \label{img:Comp1}
\end{figure}

Figure \ref{img:Comp1} shows the comparison for four different initial values.
The solid red line indicates the value of the ``worst case'' estimate,
the dash-dotted blue line our new eigenvalue estimate and the dashed
orange line the value of the finite dimensional eigenvalue $\lambda_n$.
The dotted green line indicates the ``number of modes needed'' for our
eigenvalue estimate to be valid. 
Please consider the difference between $n$, the number used in 
Theorem \ref{thm:evbound}, and $N$ the number of Fourier 
modes used for a simulation.
(e.g the condition $n\geq \sqrt{2}C_\varphi$ where $2\sqrt{2}C_\varphi + 1
= \#\text{Fourier modes needed}$).

The first two images (a) and (b) show for both our methods easy to handle
initial values, whereas (c) and (d) are only treatable with the new
eigenvalue estimate. The reason can be seen in the magnitude of the
``worst case'' estimate which amounts to around 800 in the latter examples,
whereas the new estimate stays below 200. Recall that these values are
an exponential growth-rate in our ODEs. Therefore,
an improvement of about 600 is a huge benefit.

Although it is a major improvement, this new estimate does not
resolve the problem connected to higher frequencies in the initial value
for the rigorous a-posteriori method.
This is not a huge surprise as it does not remove the exponential
growth of the error itself, it just significantly reduces its exponent.

\begin{figure}[htb]%
    \centering
    \hspace*{\fill} %
    \subfloat[$\sin(7x)$]
    {\includegraphics[width=0.45\textwidth]
        {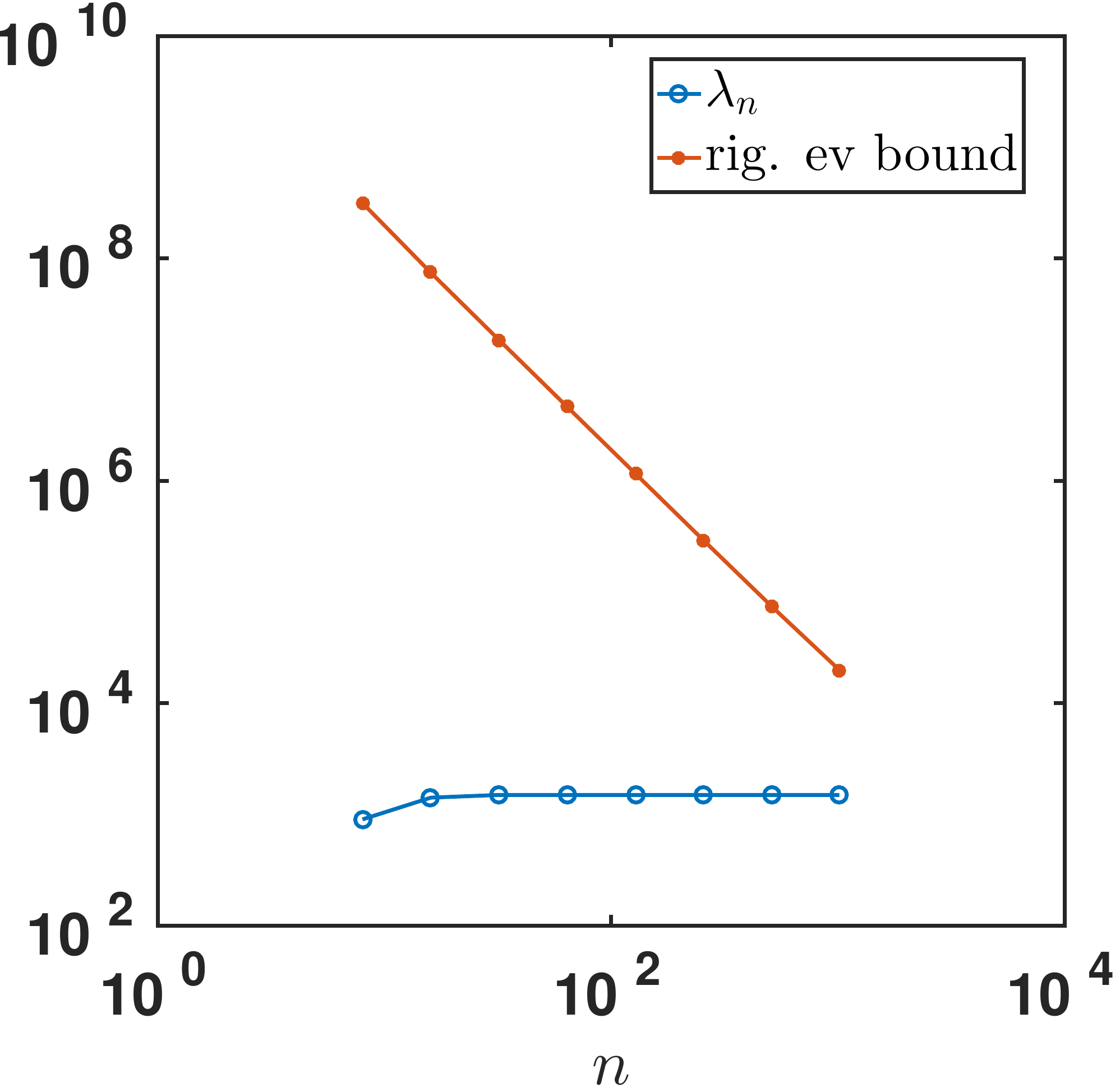}}
    \hspace*{\fill} %
    \subfloat[$\cos(2x) + \sin(2x)$]
    {\includegraphics[width=0.45\textwidth]
        {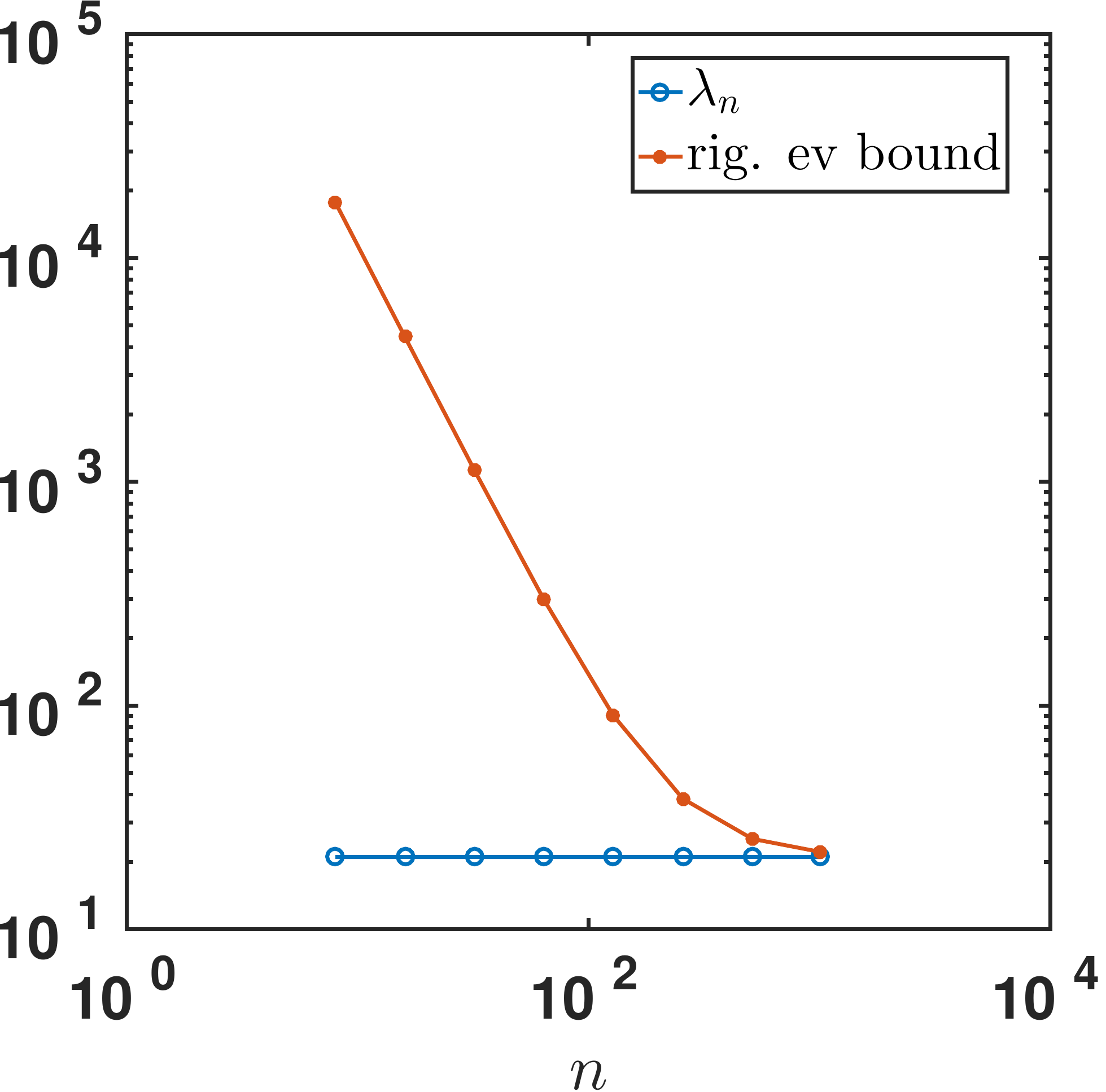}}
    \hspace*{\fill} %
    \\
    \caption
    {Convergence of the rigorous eigenvalue bound to the finite dimensional 
    eigenvalue for increasing $n$. The values for $n$ are 
    $8, 16, 32, 64, 128, 256, 512, 1024$. Please note the logarithmic 
    scale of the x- and y-axis.}
    \label{img:Convergence}
\end{figure}

In Figure \ref{img:Convergence} we can see how the rigorous eigenvalue bound 
from Theorem \ref{thm:evbound} converges to the finite dimensional 
eigenvalue $\lambda_n$ for increasing $n$. Note, that the axes are 
using a logarithmic scaling. The results show, that there is 
room for improvement if one is willing and able to use more modes 
in the eigenvalue estimate which on the other hand increases calculation time 
drastically. Also, the finite dimensional numerical eigenvalue 
stays basically constant after a certain number of modes is reached 
(i.e. that $\varphi^2$ can be represented).

\begin{figure}[htb]%
    \centering
    \hspace*{\fill} %
    \subfloat[Smallness Method 1]
    {\includegraphics[width=0.3\textwidth]
        {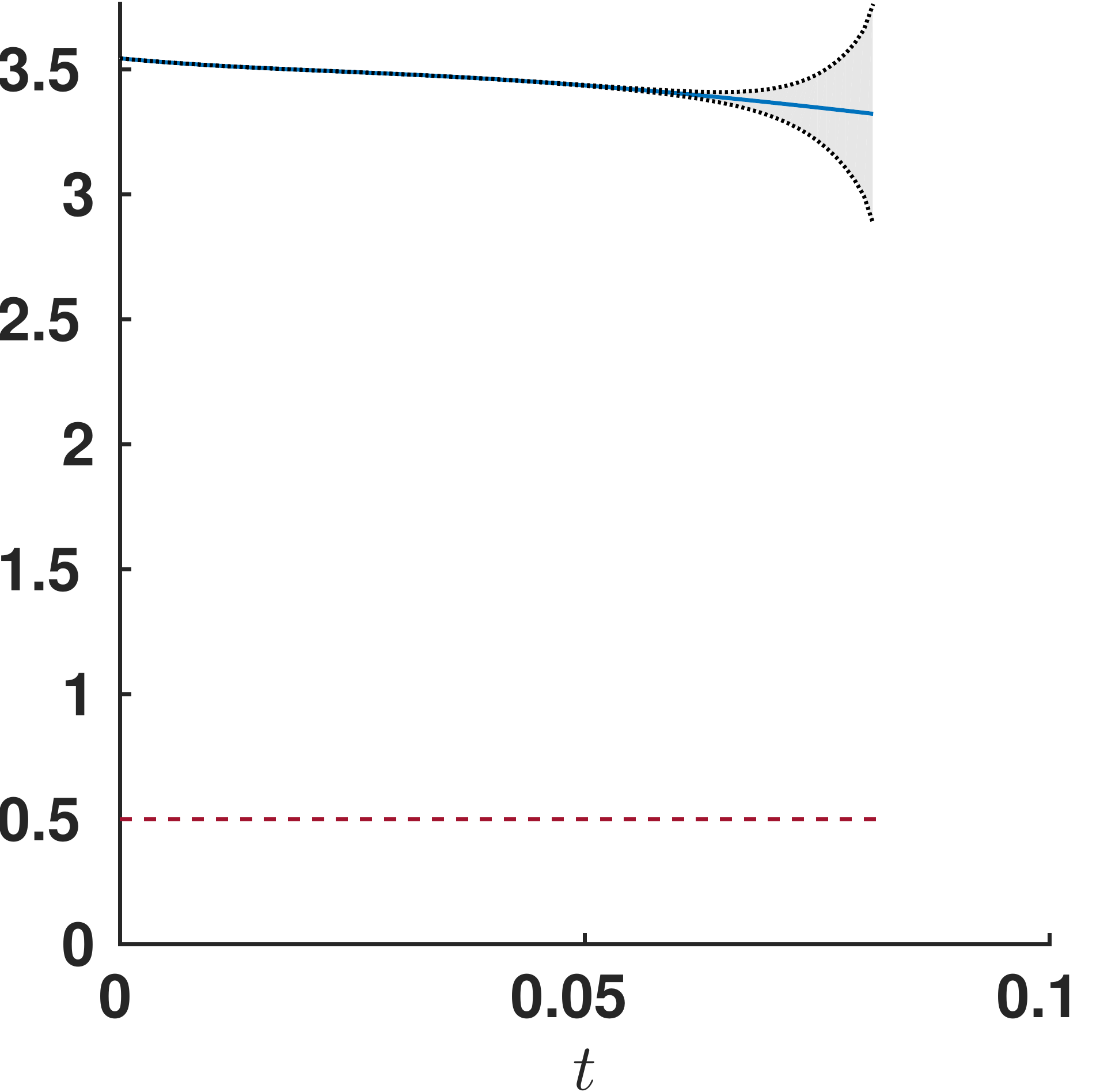}}
    \hspace*{\fill} %
    \subfloat[Smallness Method 2]
    {\includegraphics[width=0.3\textwidth]
        {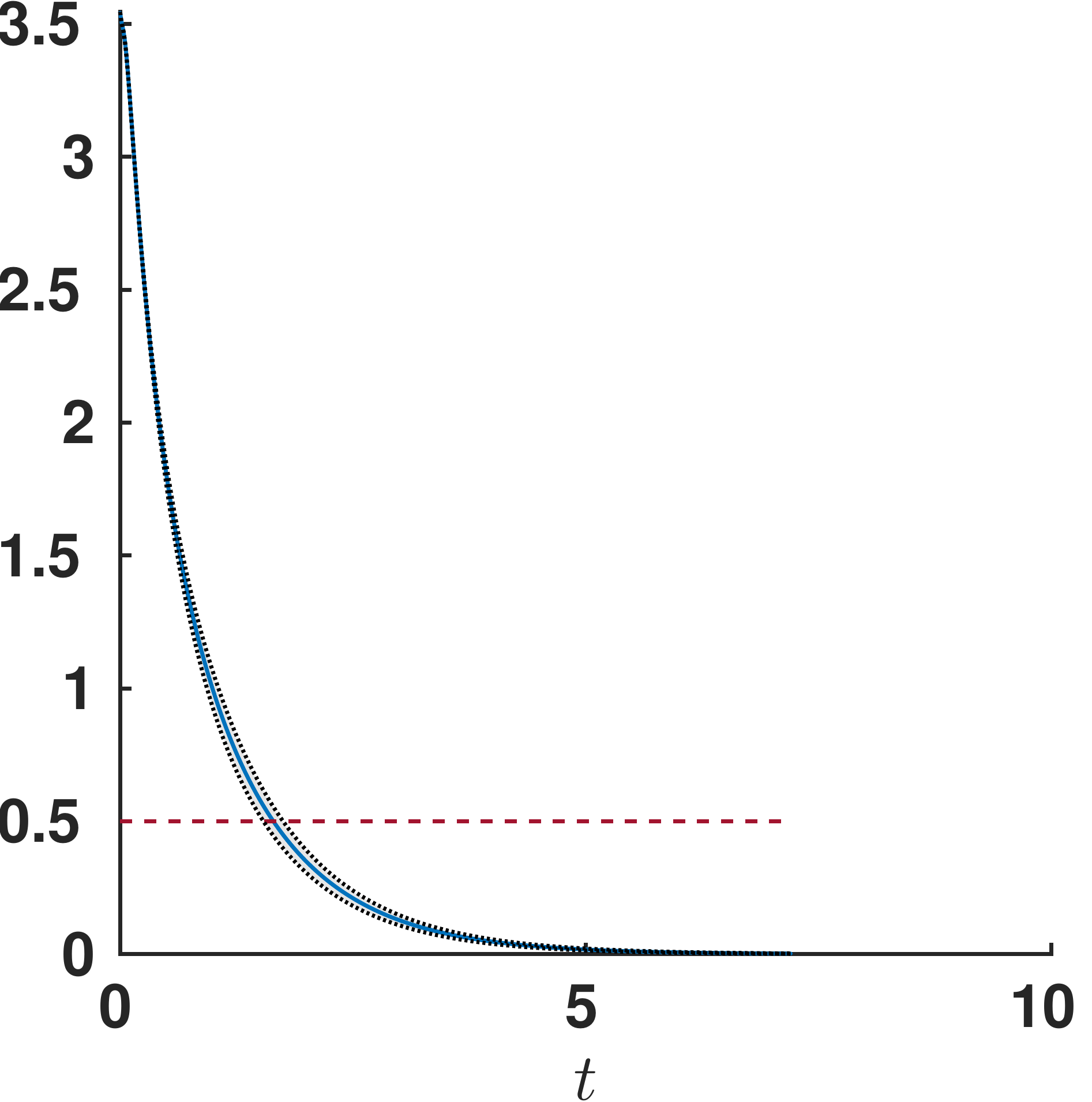}}
    \hspace*{\fill} %
    \subfloat[$\varphi$]
    {\includegraphics[width=0.3\textwidth]
        {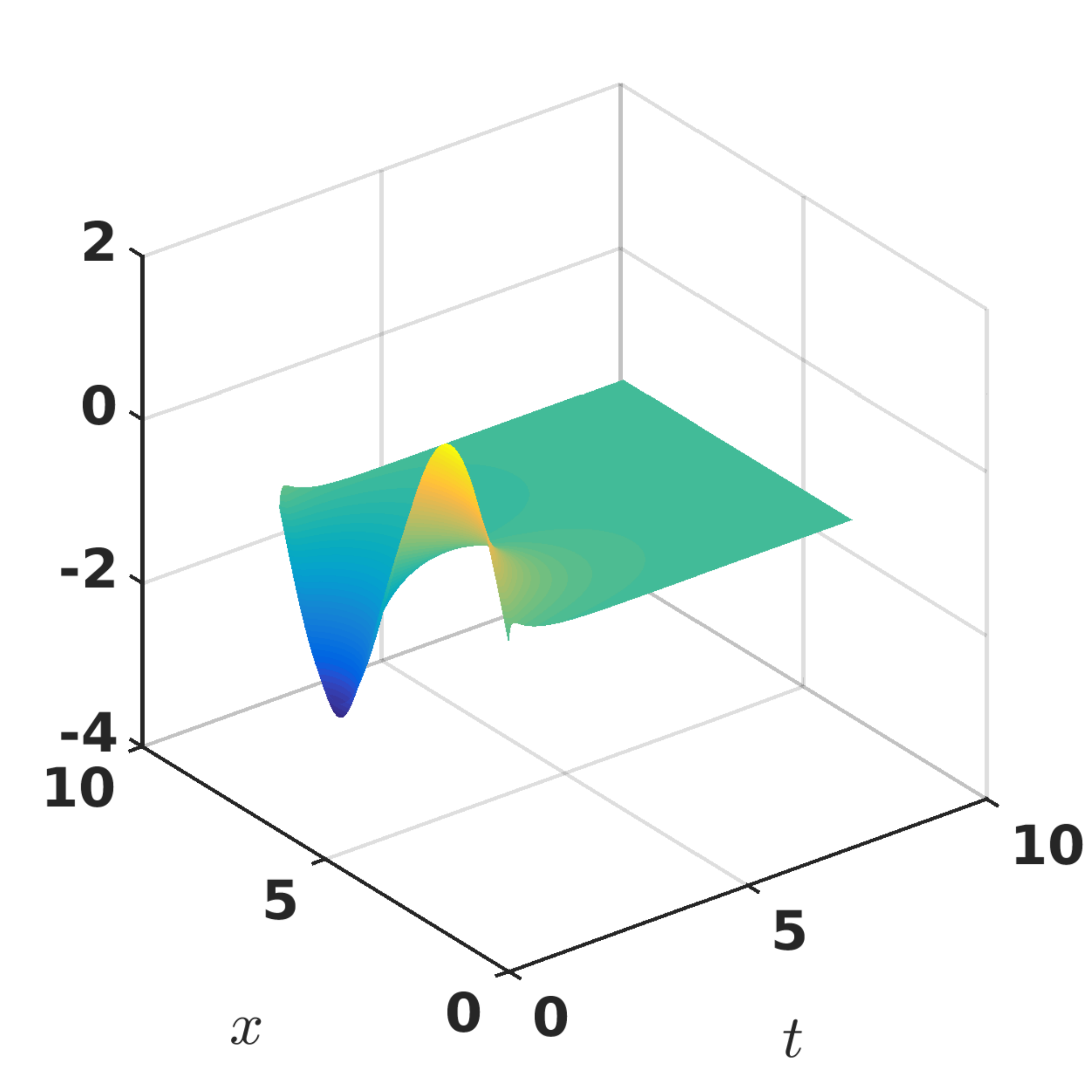}}
    \hspace*{\fill} %
    \\
    \caption[
    $u_0 = 2 \sin(x)$, 
    $N = 256$ and $h=10^{-6}$.]
    {$u_0 = 2 \sin(x)$, 
        $N = 256$ and $h=10^{-6}$. $N$ is larger than the maximum 
        for modes needed, 
        so that Method 2 (with eigenvalue estimate) is valid. 
        Method 1 (without eigenvalue estimate) fails relatively fast 
        whereas Method 2 
        succeeds in both the smallness and time criterion.
    }
    \label{img:fullcomp}
\end{figure}

Finally, in Figure \ref{img:fullcomp} we show our methods as described
above, where Method 1 uses the former ``worst case'' estimate and 
Method 2 the new eigenvalue estimate from Theorem \ref{thm:evbound}. 
The "Smallness Method X" plots will show 
the $\mathcal{H}^1$-norm of the approximation $\varphi$ 
surrounded by the gray area 
in which the smooth solution lies (the borders are given by the respective 
method). The red dotted line in these 
plots represents the threshold for the smallness criterion.
If the upper bound of the gray area falls below this threshold, we have 
global regularity.
The simulations show that whereas Method 1 reaches a blowup relatively 
fast, Method 2 stays small 
enough to reach both, the smallness and the time criterion, due to the new
eigenvalue estimate.
The corresponding plot of the eigenvalue estimate can be found in 
Figure \ref{img:evest2sin} (truncated in time, 
but the interesting part is there).
%
\section{Conclusion}
We presented a rigorous eigenvalue estimate based on 
numerical calculations to improve our previous estimates 
which relied 
heavily on general interpolation inequalities  
for numerical verification of global uniqueness for solutions of the 
surface growth equation. Our simulations show 
that this eigenvalue estimate is a huge improvement 
and suggest that the eigenvalue bound converges to the true 
eigenvalue for $n \to \infty$. 
Please keep in mind that in order to speed up the calculations 
our simulations are not fully rigorous 
as interval arithmetic was not used, although every mathematical 
preparation was carried out. 
We only wanted to establish a proof of concept that the methods do work.
\FloatBarrier
%
\section*{Acknowledgment}
This project was supported by the 
``Deutsche Forschungsgemeinschaft'' (DFG) as part of the project
BL 535/10-1 
``Numerische A-posteriori Regularität für Lösungen eines
 Oberflächenwachstumsmodells''.
\printbibliography
%
\end{document}